\documentclass[12pt]{amsart}
\usepackage[mathscr]{euscript}
\usepackage{amsthm}
\usepackage{amssymb}
\usepackage{amsmath} 
\usepackage{verbatim}
\usepackage{amscd} 
\usepackage{  latexsym }
 \newcommand{\A}{\mathcal A}  
\swapnumbers
\newtheorem{theorem}{Theorem}[section] 
\newtheorem{lemma}{Lemma}[section] 
 \theoremstyle{definition}
\newtheorem*{remarks}{Remarks}
\newtheorem*{examples}{Examples}

\begin{document} 
\date{\today} 
\title{The natural density of some sets of   square-free numbers}  

\author{Ron Brown}
\address{Department of Mathematics\\University of Hawaii\\2565 McCarthy Mall\\Honolulu, Hawaii 96822}
\email{ron@math.hawaii.edu}
\begin{abstract}
Let $P$ and $T$ be disjoint sets of prime numbers with $T$ finite.  A simple formula is given for the natural density of the set of square-free 
numbers which are divisible by all of the primes in $T$ and by none of the primes in $P$.  If $P$ is the set of primes congruent to $r$ modulo $m$
(where $m$ and $r$ are relatively prime numbers), then this natural density is shown to be $0$.
\end{abstract}
  
\maketitle

\section{Main results}\label{s:main}

In $1885$ Gegenbauer proved that the natural density of the set of square-free integers, i.e.,  the proportion of natural numbers which are square-free, is   $6/\pi^2$ \cite[Theorem 333;  reference on 
 page 272]{HW}.    In 2008   J. A. Scott conjectured  
that the proportion of natural numbers which are  odd square-free numbers is $4/\pi^2$ or, equivalently, the proportion of natural numbers which are square-free and  divisible by $2$ is $2/\pi^2$ \cite{Scott}.  The conjecture was proven in 2010 by G. J. O. Jameson, in an argument adapted from one computing
 the natural density  of the set of all square-free numbers \cite{Jameson}.    In this note we use the classical result for all square-free numbers to reprove Jameson's result and indeed to generalize it:

\begin{theorem}\label{main}
Let $P$ and $T$ be disjoint sets of prime numbers with $T$ finite.  Then the proportion of all  numbers which are square-free and divisible by 
all the primes in $T$ and by none of the primes in $P$ is 
\[
\frac{6}{\pi^2}\prod_{p \in T}\frac{1}{1+p }\prod_{p \in P}\frac{p}{1+p}.
\]
\end{theorem} 

As in the above theorem, throughout this paper   $P$ and $T$ will be disjoint sets of prime numbers with $T$ finite. The letter  $p$ will always denote a prime number.
The  term \textit{numbers} will always refer to positive integers.  Empty products, such as occurs in the first product above when $T$ is empty, are understood to equal $1$. 
If $P$ is infinite, we will argue below that the second product above is well-defined.

\begin{examples}
1.   Setting $P = \{2\}$ and $T$ equal to the empty set in the theorem we see that the natural density of the set of odd square-free numbers is $\frac{6}{\pi^2}\frac{2}{2+1} = \frac{4}{\pi^2}$; taking $T = \{2\}$ and $P$ equal to the empty set we see that the natural density of the set of even square-free numbers is  $\frac{6}{\pi^2}\frac{1}{2+1} = \frac{2}{\pi^2}$.  Thus one third of the square-free numbers  are even and two thirds are odd.  (These are Jameson's results of course.) 

2. Set  $T = \{2,3,5\}$ and $P=\{7\}$ in the theorem.  Then the theorem says that the natural density of the set of square-free numbers divisible by 30 but not by 7 is $\frac{6}{\pi^2}\frac{1}{2+1}\frac{1}{3+1}\frac{1}{5+1}\frac{7}{7+1}$, so the proportion of square-free numbers which are divisible by $30$ but not by $7$ is 
$\frac{1}{2+1}\frac{1}{3+1}\frac{1}{5+1}\frac{7}{7+1}= \frac{7}{576} $. 
\end{examples}

Our interest in the case that $P$ is infinite arose in part from a question posed by Ed Bertram: what is the natural density of the set of square-free numbers none of which is divisible by a prime congruent to $1$ modulo $4$?  The answer is zero; more generally  we have the:

\begin{theorem}\label{arithprog}
Let $r$ and $m$ be relatively prime numbers.  Then the natural density of the set of square-free numbers divisible by no prime congruent to $r$ modulo $m$ is zero. 
\end{theorem}

This theorem is a corollary of the previous theorem since, as we shall see in Section \ref{s:arithprog},   for any $r$ and $m$ as above, 
\[
\prod_{p \equiv r\mod{m}} \frac{p}{1+p}= 0.
\]

\section{A basic lemma}\label{s:bijection}

For any real number $x$ and set $B$ of numbers, we let $B[x]$ denote the number of elements $t$ of $B$ with $t \le x$.  Recall that if $\lim_{x \to \infty} B[x]/x$ exits, then it is by definition the \textit{natural density} of
$B$   \cite[Definition 11.1]{NZ}.
 
Let $\A$ denote the set of square-free numbers.   Then we let $\A(T,P)$ denote the set of elements of $\A$ which are divisble by all elements of $T$ and by no element of $P$ (so, for example, $\A=\A(\emptyset,\emptyset)$).  The set of square-free numbers analyzed in Theorem \ref{main} is $A(T,P)$.

The next lemma  shows how   the calculation of the natural density of the sets $A(T,P)$ reduces to the calculation of the natural density of   sets of the form $A(\emptyset,S)$ and, when $P$ is finite, also  reduces to the the calculation of the natural density of sets of the form $\A(S,\emptyset)$.  

\begin{lemma}\label{bijection}
For any finite set of primes $ S $ disjoint from $T$ and from $P$   and for any real number  $x$, we have $\A(T,S\cup P)[x] = \A(T \cup S,P)[xs]$ 
where $s = \prod_{p\in S}p$. Moreover, the set $\A(T\cup S,  P)$ has a natural density if and only if $\A(T  ,P\cup S)$ has a natural density, and if $D$ is the natural density of $\A(T\cup S,  P)$, then the natural density of  $\A(T  ,P\cup S)$ is $sD$.
\end{lemma}

\begin{proof}
The first assertion is  immediate from the fact that multiplication by $s$ gives a bijection from the set of   elements of $\A(T,S\cup P)$ less than or equal to $x$ to the set of elements of $\A(T \cup S,P)$ less than or equal to $xs$. This implies that
\[
\frac{\A(T,P\cup S)[x]}{x} = s\frac{\A(T \cup S,P)[xs]}{xs}.
\]
The lemma follows by  taking the limit as $x$ (and hence $xs$) goes to infinity.
 
\end{proof}

\begin{remarks}
1. We might note that if we assume that for all $T$   the sets $A(T,\emptyset)$ have   natural densities, then it is easy to compute these natural densities.  After all, for any $T$ the set 
$\A $ is the disjoint union over all subsets $S$ of $T$ of the sets $\A(T\setminus S,S)$, so by Lemma \ref{bijection} for any real number $x$  
\[
\A[x] = \sum_{S\subseteq T} \A(T\setminus S,S)[x]  = \sum_{S\subseteq T}\A(T,\emptyset)[d_Sx]
\] 
where for any $S\subseteq T$ we set $d_S = \prod_{p \in S} p$. Hence
\[
\frac{\A[x]}{x} = \sum_{S\subseteq T} d_S\frac{\A(T,\emptyset)[d_Sx]}{d_Sx}.
\]
Taking the limit as $x$ goes to $\infty $ we have 
$
6/\pi^2 = (\sum_{S\subseteq T} d_S) G
$
where $G$ denotes the natural density of $A(T,\emptyset)$.   But
\[
\sum_{S\subseteq T} d_S = \sum_{d|d_T} d = \prod_{p \in T} (1+p)
\]
\cite[Theorem 4.5]{NZ}, so indeed 

\[
G= \frac{6}{\pi^2} \prod_{p \in P} \frac{1}{ 1+p}.
\]

2. In the language of probability theory, Theorem \ref{main} says that the probability that a number in $\A$ is divisible by a prime $p$ is $1/(p+1)$ (so the probability that it is not is $p/(p+1)$) and, moreover, for any finite set $S$ of primes not equal to $p$, being divisible by $p$ is independent of being divisible by all of the elements of  $S$.
 \end{remarks}

\section{Proof of Theorem \ref{main}     when
  $P$ is finite}\label{mainproof1}

First suppose that $P = \emptyset$. We will prove Theorem \ref{main} in this case by induction on the number of elements of $T$.  The next lemma gives the induction step.

\begin{lemma}
Let $p$ be a prime number  not in $T$.  If the set $\A(T,\emptyset)$ has natural density $D$, then the set $\A(\{p\}\cup T,\emptyset) $ has natural density $ \frac{1}{p+1}D$.
\end{lemma}

\begin{proof}
For any real number $x$ we set  $E(x)  = \A(\{p\}\cup T,\emptyset)[x]$.  Let $\epsilon > 0$. The theorem says that 
\[
\lim_{x \to \infty}\frac{E(x)}{x} =  \frac{1}{p+1} D. 
\] 
Therefore it suffices to show for all choices of  $\epsilon$   above that,  for all sufficiently large $x$ (depending on $\epsilon$),
\[ 
\left| \frac{E(x)}{x} -  \frac{1}{p+1} D\right| < \epsilon.
\]

Note that $\A( T,\emptyset)$ is the disjoint union of $\A(\{p\}\cup T,\emptyset)$ and $\A(T,\{p\})$.  Hence by Lemma \ref{bijection} (applied to $\A(T,\{p\})$) for any real number $x$, 
\[
A(T,\emptyset)[x/p] = E(x/p) + E(x)
\]
and so  by the choice of $D$ there exists a number $M$ such that if $x>M$ then 
\[ 
\left| \frac{E(x)}{x/p} + \frac{E(x/p)}{x/p} - D  \right| < \epsilon/3.
\] 

We next pick an even integer $k$ such that $\frac{1} {p^k} < \frac{\epsilon}{3}$. Then
\begin{equation}\label{e2}
\left|   E(x/p^k) \right| \le     x/p^k     < \frac{\epsilon}{3} x
\end{equation}
and also (using the usual formula for summing a geometric series)
\begin{equation}\label{e1}
\left| -Dx \sum_{i=1}^{k} (-\frac{1}{p})^i - Dx\frac{1}{p+1} 
   \right| 
 =  Dx \left|\frac{(-\frac{1}{p})-(-\frac{1}{p})^{k+1}}{1-(-\frac{1}{p})} +\frac{1}{p+1}\right|  
\end{equation}
\[
 =Dx
 \left|\frac{-1+\frac{1}{p^k} }{p+1} +\frac{1}{p+1}\right|< \frac{1}{p^k} Dx < \frac{\epsilon}{3}\frac{6}{\pi^2} x  <\frac{\epsilon}{3} x  .
 \]

Now suppose that $x > p^kM$.  Then for all $i \le k$ we have $x/p^i > M$ and hence (applying the choice of $M$ above),
\begin{equation}\label{e3} 
\left| E(x) + E(x/p) - D\frac{x}{p}\right| < \frac{\epsilon}{3} \frac{x}{p}
\end{equation} 
and similarly 
\[ 
\left|- E(x/p) -E(x/p^2) + D\frac{x}{p^2 }\right| < \frac{\epsilon}{3} \frac{x}{p^2}
\]
and   
\[ 
\left| E(x/p^2) + E(x/p^3) - D\frac{x}{p^3}\right| < \frac{\epsilon}{3} \frac{x}{p^3}
\]
and 
 \[  
\left|- E(x/p^3) -E(x/p^4) + D\frac{x}{p^4}\right| < \frac{\epsilon}{3} \frac{x}{p^4}
\]

 \centering{\LARGE{$\vdots$}}

 \enlargethispage{\baselineskip}
\begin{equation} \label{ek}
\left|- E(x/p^{k-1}) -E(x/p^k) +  D\frac{x}{p^k} \right| < \frac{\epsilon}{3} \frac{x}{p^k}.
 \end{equation}
\begin{flushleft}
Using the triangle inequality to combine the inequalities (\ref{e2}) and  (\ref{e1}) together with all those between (\ref{e3}) and 
(\ref{ek}) (inclusive) and  dividing through by $x$, we can conclude that \end{flushleft} 
\[
\left| \frac{E(x)}{x} - \frac{1}{p+1} D\right| < \frac{\epsilon}{3}\left( \sum_{i=1}^{k} \frac{1}{p^i}\right) + \frac{\epsilon}{3} +  \frac{\epsilon}{3}<\epsilon.
\]
\end{proof}

Theorem \ref{main} now follows in the case that $P$ is empty from  the above lemma by induction on the number of elements of $T$.  That it is true when 
 $P$ is finite  but not necessarily empty  follows from Lemma  \ref{bijection}: in the statement of the lemma replace $P$ by $\emptyset$ and $S$ by $P$; then the natural density of $\A(T,P)$ is 
\[
\frac{6}{\pi^2}\prod_{p \in P}p   \prod_{p \in T \cup P} \frac{1}{1+p} =  \frac{6}{\pi^2}\prod_{p \in T}\frac{1}{1+p}\prod_{p \in    P}\frac{p}{1+p}.
\]

\section{Proof of Theorem \ref{main} when $P$ is infinite}\label{s:mainproof}  
 
We begin by proving the theorem in the case that $T$ is empty.  Let $p_1,p_2,p_3, \cdots$ be the strictly increasing sequence of elements of $P$. Since all the quotients $p_i/(1+p_i)$ are less than $1$, the partial products of the infinite product $\prod_i p_i/(1+p_i)$ form a strictly decreasing sequence bounded below by $0$; thus $\prod_{p \in P} p/(1+p) $ converges and its limit, say $\alpha$, is independent of the order of the factors.

First suppose that $\alpha \ne 0$. Then $\sum_{p\in P} 1/p < \infty$.  After all, we have
\[
- \log \alpha = - \sum_{p \in P} \log \frac{p}{1+p} = \sum_{p \in P} \log (1+p) -\log p > \sum_{p\in P}\frac{1}{1+p} > \frac{1}{2}\sum_{p\in P} 1/p.
\]

Now observe that $\A \setminus \A(\emptyset , P )$ is the disjoint union 
\[
 \A \setminus \A(\emptyset , P ) = \cup_{k \ge 1} \A(\{p_k\}, \{ p_1, \cdots, p_{k-1}\})
\]
since for all $b \in \A \setminus \A(\emptyset , P )$ there exists a least $k$ with $p_k|b$, so that $b\in \A(\{p_k\}, \{ p_1, \cdots, p_{k-1}\})$.

For all $n$ and $k$ we have
\[
\frac{
\A(\{p_k\}, \{ p_1, \cdots, p_{k-1}\})[n]}{n} \le \frac{|\{j: 1 \le j \le n,  p_k | j\}|}{n} \le \frac{1}{p_k}.
\]
Hence by Tannery's theorem (see \cite[p. 292]{tannery} or \cite[p. 199]{hofbauer})
the natural density of $\A \setminus \A(\emptyset,P)$ is
\[
\lim_{n \to \infty} \frac{(\A \setminus \A(\emptyset,P))[n]}{n} = \lim_{n \to \infty} \sum_{k=1}^{\infty}\frac{\A(\{p_k\}, \{ p_1, \cdots, p_{k-1}\})[n]}{n}
\]
\[
= \sum_{k=1}^{\infty}\lim_{n \to \infty}\frac{\A(\{p_k\}, \{ p_1, \cdots, p_{k-1}\})[n]}{n}
=\sum_{k=1}^{\infty}\frac{6}{\pi^2} \frac{1}{1+p_k}\prod_{i<k}\frac{p_i}{1+p_i}
\]
by the proof in the previous section of the theorem in the case that $P$ is finite.  
Writing $1/(1+p_k) = 1 - p_k/(1+p_k)$ we can see that the natural density of $\A \setminus \A(\emptyset,P)$ is therefore a limit of a telescoping sum 
\[
\frac{6}{\pi^2}\lim_{L \to \infty} \sum_{k=1}^{L} \left(\prod_{i<k}\frac{p_i}{1+p_i} - \prod_{i<k+1}\frac{p_i}{1+p_i}\right)
\]
\[
 =\frac{6}{\pi^2}\lim_{L \to \infty}\left(1- \prod _{i \le L}\frac{p_i}{1+p_i}\right) =\frac{6}{\pi^2}\left(1-\prod_{p\in P}\frac{p}{1+p}\right)
\] 
and therefore
the natural density of $\A(\emptyset,P)$ is 
 \[
\frac{6}{\pi^2} - \frac{6}{\pi^2}\left(1-\prod_{p\in P}\frac{p}{1+p}\right)=\prod_{p\in P}\frac{p}{1+p},
\]
 as was claimed.

We now consider the case that  $\alpha =\prod_{p \in P} p/(1+p)= 0$ (still assuming that $T = \emptyset$).  Suppose that $\epsilon >0$.  By hypothesis there exists a number $M$ with $\frac{6}{\pi^2} \prod_{i\le M}\frac{p_i}{1+p_i} < \epsilon/2$. Then by our proof of the theorem in the case that $P$ is finite there exists a number $L$ such that if $n>L$ then 
\[ 
\frac{(\A(\emptyset, \{p_1,p_2, \cdots, p_M\})[n]}{n} < \frac{\epsilon}{2} +\frac{6}{\pi^2}\prod_{i\le M}\frac{p_i}{1+p_i} .
\]
Then if $n>L$ we have
\[
0 < \frac{\A(\emptyset, P)[n]}{n} \le \frac{\A(\emptyset, \{ p_1,p_2, \cdots, p_M \})[n]}{n} < \frac{\epsilon}{2}+\frac{\epsilon}{2}
= \epsilon.
\]
Therefore, 
\[ 
\lim_{n \to \infty}\frac{\A(\emptyset,P)[n]}{n} = 0 = \prod_{p\in P}\frac{p}{1+p}.
\]

This completes the proof of the theorem in the case that $T= \emptyset$.  
The general case where $T$ is arbitrary then follows from Lemma \ref{bijection}, applied with $T$ and $S$ replaced respectively by $\emptyset$ and $T$: then if we set  $s = \prod_{p \in T}p$ we have that the natural density of $\A(T,P)$ equals $1/s$ times the natural density of $\A(\emptyset, P\cup T)$, i.e.,
equals
\[
\frac{6}{\pi^2}\prod_{p\in T}\frac{1}{p} \prod_{p \in T}\frac{p}{1+p}\prod_{p \in P}\frac{p}{1+p} =\frac{6}{\pi^2} \prod_{p \in T}\frac{1}{1+p}\prod_{p \in P}\frac{p}{1+p},
\]
which completes the proof of Theorem \ref{main}.

\section{Proof of Theorem \ref{arithprog}}\label{s:arithprog}

It suffices by Theorem \ref{main} to prove that $\prod_{p \in P} p/(1+p) = 0$.  Let $M$ be any number. By a lemma of K. K.  Norton
\cite[Lemma 6.3]{norton} there exists a constant $B$ independent of the choices of $r,m,$ and $M$ such that 
\[
\left|\sum_{M>p\in P}\frac{1}{p} - \frac{\log\log M}{\phi(m)} \right| < B \frac{\log(3m)}{\phi(m)}.
\]
As we observed earlier, for any $p$ we have 
\[ 
1/(2p) < 1/(1+p) < \log(1+p) - \log p
\]
 so
\[
\log \prod_{M > p \in P} \frac{p}{1+p} = \sum_{M>p \in P}(\log p -\log(1+p)) 
\]
\[
< - \sum_{M>p\in P}\frac{1}{2p} < \frac{1}{2}\left(-\frac{\log\log M}{\phi (m)} + B \frac{\log(3m)}{\phi(m)}\right)
\]
which has  limit  $-\infty$ as $M \to \infty$.  Hence 
\[
\prod_{p \in P} \frac{p}{1+p}=
\lim_{M \to \infty} \prod_{M>p \in P}\frac{p}{1+p} = 0,
\]
completing the proof of Theorem \ref{arithprog}.

\bibliographystyle{plain}
\bibliography{\jobname} 

\begin{thebibliography}{00}  


 \bibitem{HW} 
  G.H. Hardy and E.M. Wright, \textit{An Introduction to the Theory of Numbers (5th ed.)}  (Oxford
Univ. Press, 1979).

 \bibitem{NZ} I. Niven and H. S. Zuckerman, \textit{An Introduction to the Theory of Numbers (4th ed.)} (Wiley, 1980).

  
\bibitem{Scott}  J.A. Scott, Square-free integers once again, \textit{Math. Gazette} {\bf 92} (2008), 70 -- 71.

\bibitem{Jameson} 
G. J. O. Jameson, 
 Even and odd square-free numbers, \textit{Math. Gazette} {\bf 94} (2010), 123--127. 
 
 \bibitem{norton}
K. K. Norton, On the number of restricted prime factors of an integer, \textit{Illinois J. Math.} {\bf 20} (1976), 681 -- 705. 

\bibitem{tannery}
 J. Tannery, Introduction \`a la Th\'eorie des Fonctions d'une Variable, 2 ed., \textit{Tome 1, Libraire Scientifique A} (Hermann, Paris, 1904).
 
 \bibitem{hofbauer}
J. Hofbauer,  A simple proof of $1 + 1/22 + 1/32 + \cdots = \pi^2/6$ and related identities, \textit{ American Math. Monthly} {\bf 109} (2002), 196 -- 200.
  


\end{thebibliography}
 
\end{document}